\documentclass[10pt, twoside, a4paper, fleqn]{article}

\usepackage{latexsym}
\usepackage{amscd}
\usepackage{theorem}
\usepackage{pifont}
\usepackage{amsfonts}
\usepackage{xspace}
\usepackage{amssymb}
\usepackage{fancyhdr}
\usepackage{mathrsfs}
\usepackage{amsmath}%
\usepackage{graphics}%
\usepackage{graphicx}%
\usepackage{tikz-cd}
\usepackage{euscript}%
\usepackage{psfrag}%
\usepackage{empheq}%
\usepackage{upgreek}
\usepackage[pagewise]{lineno}
\DeclareMathAlphabet{\mathpzc}{OT1}{pzc}{m}{it}

{\theorembodyfont{\slshape} \newtheorem{theorem}{Theorem}[section]}
{\theorembodyfont{\slshape} \newtheorem{definition}{Definition}[section]}
{\theorembodyfont{\slshape} \newtheorem{lemma}[theorem]{Lemma}}
{\theorembodyfont{\slshape} }
{\theorembodyfont{\slshape} \newtheorem{corollary}[theorem]{Corollary}}
{\theorembodyfont{\slshape} \newtheorem{remark}[theorem]{Remark}}
{\theorembodyfont{\slshape} \newtheorem{example}[theorem]{Example}}

\newcommand{\reals}{\mathbb{R}}

\newcommand{\R}{\mathbb{R}^{n}_{+}}
\newcommand{\G}{\mathcal{G}}
\newcommand{\D}{D_{A}}
\newcommand{\N}{N_{A}}
\newcommand{\C}{\lbrace A_{1},\,A_{2},\,\cdots ,\,A_{N}\rbrace}
\newenvironment{proof}{\noindent\textbf{Proof:\ }}{$\hfill{\Box}$}

\numberwithin{equation}{section}

\title{\textsc{Dynamics of Products of Matrices in Max Algebra}}    

\author{Sachindranath Jayaraman\footnote{Corresponding Author} \\ {\tt sachindranathj@iisertvm.ac.in} \bigskip 
	\\ Yogesh Kumar Prajapaty \\ {\tt prajapaty0916@iisertvm.ac.in} \bigskip 
	\\ Shrihari Sridharan \\ {\tt shrihari@iisertvm.ac.in} \bigskip 
 \bigskip \\ {\sl Indian Institute of Science Education and Research}
  \\ {\sl Thiruvananthapuram (IISER-TVM), India.}} 

\date{September 28, 2021}

\begin{document}

\maketitle
\thispagestyle{empty}

\begin{abstract} 
\noindent 
The aim of this manuscript is to understand the dynamics of matrix products in a max algebra. A consequence of the Perron-Fr\"{o}benius theorem on periodic points of a nonnegative matrix is generalized to a max algebra setting. The same is then studied for a finite product associated to a $p$-lettered word on $N$ letters arising from a finite collection of nonnegative matrices, with each member having its maximum circuit geometric mean at most $1$.
\end{abstract}

\begin{tabular}{r c l} 
{\bf Keywords} & : & Products of nonnegative matrices; \\ 
& & Max-algebras; \\ 
& & Boolean matrices; \\ 
& & Fr\"{o}benius normal form of nonnegative matrices; \\ 
& & Circuit geometric mean; \\
& & Periodic points. \\ 
& & \\ 
& & \\ 
{\bf MSC 2010 Subject} & : & 15A80; 15B34; 37H12. \\ 
{\bf Classifications} & & 
\end{tabular} 
\bigskip 

\section{Introduction} 
\noindent
We work throughout over the field $\mathbb{R}$ of real numbers and use the standard notations $\mathbb{R}^{n}$ and $M_{n}(\mathbb{R})$ to denote the $n$-dimensional real vector space of $n$-tuples of real numbers and the real vector space of $n \times n$ square matrices with real entries respectively. We concern ourselves with only those matrices whose entries are nonnegative real numbers. This set will be denoted by $M_{n} (\mathbb{R}_{+})$. Other notations and terminologies used in this work will be introduced later. 
\medskip 

\noindent 
Recall that given a self map $f$ on a topological space $X$, an element $x \in X$ is called a \emph{periodic point} of $f$ if there exists a positive integer $q$ such that $f^{q} (x) = x$. In such a case, the smallest such integer $q$ that satisfies $f^{q} (x) = x$ is called the \emph{period} of the periodic point $x$. The starting point for this work is the following consequence of the Perron-Fr\"{o}benius theorem. 

\begin{theorem} (Theorem B.4.7, \cite{LN-Book}) 
\label{starting-theorem}
Let $A \in M_{n} (\mathbb{R}_{+})$ with spectral radius less than or equal to $1$. Then, there exists a positive integer $q$ such that for every $x \in \mathbb{R}^{n}$ with $\left( \Vert A^{k} x \Vert \right)_{k\, \in\, \mathbb{N}}$
 bounded, we have 
\[ \lim\limits_{k\, \to\, \infty} A^{k q} x\ \ =\ \ \xi_{x}, \] 
where $\xi_{x}$ is a periodic point of $A$ whose period divides $q$. 
\end{theorem}

\noindent 
In an attempt to generalize Theorem \eqref{starting-theorem}, when the matrix $A$ in the above theorem is replaced by a product of the matrices $A_{r}$'s, possibly an infinite one, drawn from the finite collection of nonnegative matrices, $\big\{ A_{1},\, A_{2},\, \ldots,\, A_{N} \big\}$, the following result was obtained. The details may be found in \cite{SJ-YP-SS-1}. 

\begin{theorem}\cite{SJ-YP-SS-1} 
\label{SJYPSS-1}
Let $\big\{ A_{1},\, A_{2},\, \ldots,\, A_{N} \big\},\ N < \infty$, be a collection of $n \times n$ matrices with nonnegative entries, each having spectral radius atmost $1$. Assume that the collection has a nontrivial set of common eigenvectors, $E$. For any finite $p$, let $A_{\omega}$ denote the matrix product associated to a $p$-lettered word $\omega$ on the letters $\big\{ 1, 2, \ldots, N \big\}$. Suppose $\mathcal{LC} (E)$ denotes the set of all real linear combinations of the vectors in $E$. Then, for any $x \in \mathcal{LC} (E)$, there exists an integer $q \geq 1$ (independent of $\omega$) such that 
\begin{equation} 
\label{Limit}
\lim\limits_{k\, \to\, \infty} A_{\omega}^{k q} x\ \ =\ \ \xi_{x}, 
\end{equation}
where $\xi_{x}$ (independent of $\omega$) is a periodic point of $A_{\omega}$, whose period divides $q$. 
\end{theorem}

\noindent 
The aim of this work is to explore the possibilities of extending Theorems \eqref{starting-theorem} and \eqref{SJYPSS-1} in the setting of max algebras. By a max algebra, we mean the triple $(\mathbb{R}_{+}, \oplus, \otimes)$, where $\mathbb{R}_{+}$ denotes the set of nonnegative real numbers, $\oplus$ denotes the binary operation of taking the maximum of two nonnegative numbers and $\otimes$ is the usual multiplication of two numbers. There are several abstract examples of max algebras. The one given above is more amenable to work with, while dealing with nonnegative matrices. Another example is the set of real numbers, together with $-\infty$, equipped with the binary operations of maximization and addition, respectively. The latter system is isomorphic to the former one via  the exponential map. A good reference on max algebras is the monograph by Butkovic \cite{Butkovic-Book}. For Perron-Fr\"{o}benius theorem in max algebras, one may refer to \cite{Bapat}. 

\section{Preliminaries on max algebras} 
\label{prelims} 

We introduce some preliminary notions from max algebras that we use, in this section. In particular, we introduce the notion of a max eigenvalue, the corresponding max eigenvector and a few properties of the same. We begin with a description of the matrix product in this max algebra. Recall that for $ a,b \in \reals_{+}$ 
\[ a\ \oplus\ b\ \ =\ \ \max\{ a,\ b \} \quad \text{and} \quad a\ \otimes\ b\ \ =\ \ a b. \] 
Let $A$ and $B$ be two $n \times n$ nonnegative matrices. Then the matrix product of $A$ and $B$ is defined by 
\[ [ A \otimes B ]_{ij}\ \ =\ \ \max_{k}\ \lbrace a_{ik}\ \otimes\ b_{kj} \rbrace. \] 

\noindent 
We now define the notion of reducible and irreducible matrices in usual matrix product. 

\begin{definition}
An $n \times n$ matrix $A$ is said to be reducible if there exists a permutation matrix $P$ such that 
\begin{equation} 
\label{reducible}
P A P^{T}\ \ =\ \  \left[\begin{matrix}
							B	&	C\\
							0	&	D	
						\end{matrix}\right],
\end{equation}
where $B$ and $D$ are square matrices. Otherwise $A$ is irreducible.
\end{definition}

\noindent
Note that the notion of irreducibility is the same in the max algebra and the usual Euclidean algebra. If $A$ is reducible and is in the form \eqref{reducible}, and if a diagonal block is reducible, then this block can be reduced further via permutation similarity. Continuing this process, we have a suitable permutation matrix $P$ such that $PAP^{T}$ is in the block triangular form 
\begin{equation} 
\label{Frobenius normal form}
PAP^{T}\ \ = \ \  \left[\begin{matrix}
						A_{11}	&	A_{12}	& \cdots	& A_{1m}\\
						0		&	A_{22}	& \cdots	& A_{2m}\\
						\vdots	&	\vdots	& \ddots	& \vdots\\
						0		&	0		& \cdots	& A_{mm}\\
						\end{matrix}\right],
\end{equation}
where each block $A_{ii}$ is square and is either irreducible or a $1 \times 1$ null matrix. This block triangular form is called the \emph{Fr\"{o}benius normal form}.
\medskip

\noindent
Given an $n \times n$ nonnegative matrix $A$, there is a natural way to associate a simple weighted directed graph $\G(A)$ to the matrix as follows: $\G(A)$ has $n$ vertices, say $1, \ldots, n$, such that there is an edge from $i$ to $j$ with weight $a_{ij}$ if and only if $a_{ij} > 0$. By a circuit, we always mean a simple circuit. In contrast, our paths may include a vertex and/or an edge more than once. By the product of a path, we mean the product of the weights of the edges in the path.
\medskip

\noindent
Let $(i_{1}, i_{2}), (i_{2}, i_{3}), \ldots, (i_{k}, i_{1})$ be a circuit in $\G(A)$. Then $a_{i_{1} i_{2}} a_{i_{2} i_{3}} \cdots a_{i_{k} i_{1}}$ is the corresponding circuit product and its $k^{th}$ root is the circuit geometric mean corresponding to the circuit. The maximum among all possible circuit geometric means in $\G(A)$ is denoted by $\mu(A)$. A circuit with circuit geometric mean equal to $\mu(\G(A))$ is called a \emph{critical circuit}, and vertices on the critical circuits are \emph{critical vertices}. Assuming that simultaneous row and column permutations have been performed from the principal submatrix of $A$, the critical matrix of $A$, denoted by $A^{C}=[a_{ij}^{C}],$ is formed from the principal submatrix of $A$ on the rows and columns corresponding to critical vertices, by setting $a_{ij}^{C}=a_{ij}$ if $(i,\, j)$ is in a critical circuit, and $a_{ij} = 0$ otherwise. Thus the critical graph $\G(A^{C})$ has vertex set $V_{C}=\lbrace \text{critical vertices}\rbrace$.  
\medskip 

\noindent 
Note that the weighted directed graph associated with $A$, namely $\mathcal{G}(A)$ and the weighted directed graph associated with $PAP^{T}$, for a permutation matrix $P$, namely $\mathcal{G}(PAP^{T})$ are isomorphic. This implies that the corresponding circuit geometric means are equal. For more on circuit geometric means, see \cite{Bapat, BP-book, Elsner-van}. 
\medskip

\noindent
Let $A$ be a nonnegative matrix. We say that $\lambda$ is a \emph{max eigenvalue} of $A$ if there exists a nonzero, nonnegative vector $x$ such that $A \otimes x = \lambda x$. Further, $x$ is called a \emph{max eigenvector} associated to $\lambda$. The following result is referred to as the max version of the Perron-Fr\"{o}benius theorem and can be found in \cite{Bapat}.

\begin{theorem} (Theorem 2, \cite{Bapat}) 
\label{PF theorem maxalgebra}
Let $A$ be an $n \times n$ nonnegative, irreducible matrix. Then there exists a positive vector $x$ such that $A \otimes x = \mu(A) x$. 
\end{theorem}

\noindent 
It can be proved that an irreducible matrix has only one eigenvalue with a corresponding eigenvector in max algebras. Moreover, unlike the eigenvalues in the Euclidean algebra, the max eigenvalue has an explicit formula given by,
 \[ \mu(A)\ \ =\ \ \mu (A^{C})\ \ =\ \ \max\limits_{\ell} \left( a_{i_{1} i_{2}} a_{i_{2} i_{3}} \cdots a_{i_{\ell} i_{1}} \right)^{\frac{1}{\ell}}. \]

\noindent
We now define a relation between two vertices in $\G(A)$. 

\begin{definition}
Let A be a nonnegative square matrix of order $n$. For $1 \leq i,\, j \leq n$, we say that $i$ has access to $j$ if there is a path from vertex $i$ to vertex $j$ in $\G(A)$, and that $i$ and $j$ communicate if $i$ has access to $j$ and $j$ has access to $i$.   
\end{definition}

\noindent
Communication is an equivalence relation. Note that in the identity matrix, we assume that every vertex $i$ communicates with itself. Further, we assume by convention that for any matrix $A,\ A^{0}$ in the max algebra is nothing, but the identity matrix. The following result concerns the spectrum of a nonnegative matrix in max algebras.

\begin{theorem} \cite{Bapat} 
\label{Max version of Frobenius-Victory Theorem}
Let $A$ be an $n \times n$ nonnegative matrix in Fr\"{o}benius normal form,
\[ \left[\begin{matrix}
						A_{11}	&	A_{12}	& \cdots	& A_{1m}\\
						0		&	A_{22}	& \cdots	& A_{2m}\\
						\vdots	&	\vdots	& \ddots	& \vdots\\
						0		&	0		& \cdots	& A_{mm}\\
						\end{matrix}\right]. \]
Then $\lambda$ is an eigenvalue with a corresponding nonnegative eigenvector if and only if there exists a positive integer $i \in \lbrace 1, 2, \cdots, m \rbrace$ such that $\lambda = \mu(A_{ii})$ and furthermore, class $j$ does not have access to class $i$ whenever $\mu(A_{jj}) < \mu(A_{ii})$. 
\end{theorem}

\noindent 
Writing the Fr\"{o}benius form of $A$ as the sum of its diagonal blocks and the strict upper triangular block, say 
\begin{eqnarray*}
PAP^{T}	&	= 	& \left[ \begin{matrix}
						A_{11}	&	A_{12}	& \cdots	& A_{1m}\\
						0		&	A_{22}	& \cdots	& A_{2m}\\
						\vdots	&	\vdots	& \ddots	& \vdots\\
						0		&	0		& \cdots	& A_{mm}\\
						\end{matrix} \right] \\
		&	=	& \left[ \begin{matrix}
						A_{11}	&	0		& \cdots	& 0		\\
						0		&	A_{22}	& \cdots	& 0		\\
						\vdots	&	\vdots	& \ddots	& \vdots\\
						0		&	0		& \cdots	& A_{1m}\\
						\end{matrix}\right] \oplus \left[\begin{matrix}
						0		&	A_{12}	& \cdots	& A_{1m}\\
						0		&	0		& \cdots	& A_{2m}\\
						\vdots	&	\vdots	& \ddots	& \vdots\\
						0		&	0		& \cdots	& 0		\\
						\end{matrix} \right] \\			
		&	=:	&		\D \oplus \N, \\	
\end{eqnarray*}    
we observe that $\N$ does not contribute in computing the circuit geometric mean of $PAP^{T}$, and consequently has no role in $\mu(A)$. This follows from Theorem \eqref{Max version of Frobenius-Victory Theorem} and the communication relation defined above. Moreover, note that $\mu(A) = \mu(\D) = \max{\lbrace \mu(A_{11}),\, \mu(A_{22}),\, \ldots ,\, \mu(A_{mm})\rbrace}$. 
\medskip 

\noindent
We illustrate Theorem \eqref{Max version of Frobenius-Victory Theorem} below. Similar examples can be found in \cite{Bapat}.

\begin{example}
Let
\begin{eqnarray*}
A\ = \left[\begin{matrix}
							4	&	0\\
							0	&	5		
						\end{matrix}\right], & B\ =\   \left[\begin{matrix}
															4	&	2\\
															0	&	5
													\end{matrix}\right]& \text{and} \ \ \  C\ =\   \left[\begin{matrix}
							5	&	2\\
							0	&	4
						\end{matrix}\right].
\end{eqnarray*}

\noindent 
Note that $A$ has two max eigenvalues, namely, $4$ and $5$, with corresponding max eigenvectors $(1,\, 0)$ and $(0,\,1)$ respectively. Observe that $B$ has two max eigenvalues namely $4$ and $5$ with corresponding max eigenvectors $(1,\, 0)$ and $(8,\, 5)$ respectively, whereas $C$ has only one max eigenvalue $5$.
\end{example}

\noindent 
We now state two useful results due to Elsner. The first one is a {\it $DAD$-type theorem} that says that a nonnegative irreducible matrix can be bounded by a matrix of all ones via a $DAD$ transform. The second one concerns the period of a nonnegative irreducible matrix with $\mu(A) = 1$. More about this will be discussed in the section on Boolean matrices.

\begin{lemma} \cite{Elsner-van 1} 
\label{DAD} 
Let $A$ be an irreducible matrix with $\mu(A) \leq 1,\ x \in \R,\ x \neq 0$ and $z = A^{*} x$, where $A^{*} = I \oplus  A \oplus A^{2} \oplus \ldots \oplus A^{n - 1}$. Then $z \in (\R)^{\circ}$. If $D = {\rm diag} (z_{1}, z_{2}, \ldots, z_{n})$, then 
\[ \left( D^{-1}AD \right)_{ij}  = \left( D^{-1} \otimes A \otimes D \right)_{ij} \leq 1. \] 
\end{lemma}
   
\begin{theorem} \cite{Elsner-van} 
\label{Elsner-result}
Assume that $A$ is a nonnegative irreducible matrix with $\mu(A)=1$. Then, there exist $q$ and $t_{0}$ in $\mathbb{Z}_{+}$ such that for all $t\geq t_{0}$, we have 
\begin{equation}
A^{t\, +\, q} \ \ =\ \ A^{t},\ \ \ \ \text{where the powers are taken in the max algebra}.
\end{equation}
\end{theorem}

\noindent
We end this section with bounds and inequalities for $\mu(A)$.

\begin{lemma} \cite{Bapat-Stanford-van} 
\label{bound1}
Let $A$ be an $n \times n$ irreducible nonnegative matrix. 
\begin{enumerate}
\item Suppose there exist $\eta_{1} > 0$ and a vector $z^{(1)} \neq 0$ such that $A \otimes z^{(1)} \geq \eta_{1} z^{(1)}$. Then $\mu(A) \geq \eta_{1}$.
\item  Suppose there exist $\eta_{2} > 0$ and a vector $z^{(2)} \neq 0$ such that $A \otimes z^{(2)} \leq \eta_{2} z^{(2)}$. Then $\mu(A) \leq \eta_{2}$.
\end{enumerate}
\end{lemma}

\noindent 
Consequently, we have: 

\begin{corollary} \cite{Bapat-Stanford-van} 
\label{mu-bound}
Let $A$ be an $n \times n$ be nonnegative matrix. Then 
\[ \min_{i}\ \max_{j}\ a_{ij}\ \ \leq\ \ \mu(A)\ \ \leq\ \ \max_{i,\, j}\ a_{ij}. \] 
\end{corollary}

\noindent 
Let $A$ and $B$ be two commuting matrices. The following result gives us the circuit geometric mean of the addition and product of the two matrices in max algebras. 

\begin{theorem} \cite{RKatz} 
\label{commuting matrices mu inequality}
Let $A$ and $B$ commute. Then,
\begin{enumerate}
\item $\mu \left( A \otimes B \right)\ \ \leq\ \ \mu(A) \otimes \mu(B)$. 
\item $\mu \left( A \oplus B \right)\ \ \leq\ \ \mu(A) \oplus \mu(B)$. 
\end{enumerate}
Moreover, equality holds in both the above inequalities if the matrices $A$ and $B$ are irreducible.
\end{theorem}

\section{Statements of results}

We state the main results of this paper, in this section. Readers are alerted that from now onwards, all the powers and products of matrices are considered in the max algebra setting, even though we do not explicitly involve a notation. The first result includes an analogue of Theorem \eqref{starting-theorem}, in the max algebra setting. 

\noindent
\begin{theorem} 
\label{result 1}
Let $A$ be a matrix with $\mu(A) \leq 1$. Then, there exists an integer $q$ such that for every $1 \le j \le q$, we have 
\[ \lim\limits_{k\, \to\, \infty} A^{k q + j}\ \ =\ \ \widetilde{A^{(j)}}, \] 
where $\widetilde{A^{(j)}}$ is a periodic matrix, whose period divides $q$. Further, for every $x \in \R$, we have 
\[ \lim\limits_{k\, \to\, \infty} A^{k q + j} \otimes x\ \ =\ \ \xi_{x}^{(j)}, \] 
where $\xi_{x}^{(j)}$ is a periodic point of the matrix $A$, whose period divides $q$. 
\end{theorem}

\noindent
The following example illustrates that the period of the periodic matrix $\widetilde{A^{(j)}}$ and the period of the periodic point $\xi_{x}^{(j)}$ have no correlation between them. Consider $A\ = \left[ \begin{matrix} 0.5 & 1 \\ 1 & 0.5 \end{matrix} \right]$. Then, we obtain $q = 2$. Further, $\widetilde{A^{(1)}}\ =\ A$ whose period is $2$ whereas, $\widetilde{A^{(2)}} \ = \left[ \begin{matrix} 1 & 0.5 \\ 0.5 & 1 \end{matrix} \right]$ has period $1$. For $x = (1,\, 1)^{T}$, we find that the period of the periodic point $\xi_{x}^{(j)}$ is $1$, for $j = 1, 2$, while for $y = (1,\, 0)^{T}$, the period of the periodic point $\xi_{y}^{(j)}$ is $2$, for $j = 1, 2$. 
\medskip 

\noindent 
The next result deals with a finite collection of matrices, say $\C$, with each of them having $\mu(A_{i}) \le 1$. For any finite $p\in \mathbb{N}$, we denote the set of all $p$-lettered words on the first $N$ positive integers by
\[ \Sigma_{N}^{p}\ \ :=\ \ \lbrace \omega = \left(\omega_{1}\; \omega_{2}\; \cdots\; \omega_{p} \right)\ : \ 1\; \leq\; \omega_{i}\; \leq\; N \rbrace, \] 
and for any $\omega \in \Sigma_{N}^{p}$, we define the matrix product 
\begin{equation} 
\label{matrixproduct} 
A_{\omega}\ \ :=\ \ A_{\omega_{p}}\; \otimes\; A_{\omega_{p - 1}}\; \otimes\; \cdots\; \otimes\; A_{\omega_{1}}. 
\end{equation} 

\noindent 
We now state the second result of this paper, which is analogous to Theorem \eqref{SJYPSS-1}, in the max algebra setting, however under the hypothesis that the collection of matrices $\C$ commute pairwise. 

\begin{theorem} 
\label{result 2}
Let $\C$ be a finite collection of pairwise commuting nonnegative matrices with $\mu(A_{i}) \le 1$, for every $1 \le i \le N$. Suppose $\omega \in \Sigma_{N}^{p}$ is a $p$-lettered word in which the letter $i$ occurs $p_{i}$ many times. Then, there exists an integer $q$ such that for every $1 \le j \le q$, we have 
\[ \lim\limits_{k\, \to\, \infty} A_{\omega}^{k q + j}\ \ =\ \ \bigotimes\limits_{i\, =\, 1}^{N} \left( \widetilde{A^{(j)}_{i}} \right)^{p_{i}}. \] 
\end{theorem}

\noindent 
Before stating our third result, let us assume that the collection $\C$ of matrices do not necessarily commute, nevertheless, has a set of common eigenvectors, say, $E = \left\{ v_{1},\, v_{2},\, \cdots,\, v_{m} \right\}$. Suppose $v_{i} = \left( v_{i}^{(1)},\, \ldots,\, v_{i}^{(n)} \right)$. Define  
\begin{eqnarray*} 
\mathcal{LC} (E) & = & \left\{ \bigoplus_{i\, =\, 1}^{m}\ \left(\alpha_{i} \otimes v_{i}\right)\ :\ \alpha_{i} \in \reals_{+} \right\} \\ 
& = & \left\{ u = \left( u^{(1)},\, \ldots,\, u^{(n)} \right) \in \R\ :\ u^{(j)} = \max \left\{ \alpha_{1} v_{1}^{(j)},\, \cdots,\, \alpha_{m} v_{m}^{(j)} \right\};\ \alpha_{i} \in \reals_{+} \right\}. 
\end{eqnarray*} 

\noindent 
We now state our third theorem in this article, which is analogous to Theorem \eqref{SJYPSS-1}, in the max setting. 

\begin{theorem} 
\label{result 3}
Let $\C$ be a collection of nonnegative matrices with $\mu(A_{i}) \leq 1$. Assume that the collection has a non-trivial set of common eigenvectors, say $E$. Let $A_{\omega}$ be the matrix product associated with the word $\omega \in \Sigma_{N}^{p}$. Then, for any $x \in \mathcal{LC} (E)$, we have 
\[ \lim\limits_{k\, \to\, \infty} A_{\omega}^{k } \otimes x\ \ =\ \ \xi_{x}, \] 
where $\xi_{x}$ is a fixed point for every matrix in the collection $\C$. 
\end{theorem}

\noindent
Note that Theorem \eqref{result 3} is also applicable in the case of commuting matrices. In the case of commuting matrices, Theorem \eqref{result 3} however, does not give a complete picture as against Theorem \eqref{result 2}.  
 
 \section{Boolean matrices}

Recall that $A$ is a Boolean matrix if all the entries of $A$ are either $0$ or $1$. The following fact on Boolean matrices can be found in \cite{Pang-Guu, Rosenblatt}. Let $A$ be a Boolean matrix. Then there exist $c_{0}$ and $q$ such that 
\begin{equation} 
\label{B Period}
A^{c\, +\, k q}\ =\ A^{c}, \ \ \text{for all}\ \ k \in \mathbb{N},\ \ c \geq\ c_{0}\ \geq\ 1. 
\end{equation}

\noindent 
The minimal such $q$ is called the period of the matrix $A$. Notice that this definition of the period of a Boolean matrix coincides with the {\it same} notion introduced in Theorem \eqref{Elsner-result} of Section \eqref{prelims}. We mean this same notion of periodicity in Theorem \eqref{result 1}. We now introduce a few more notions that will help us in understanding more about the period of a matrix. These are taken from \cite{Rosenblatt}.

\begin{definition} 
\label{fourpointone}
A Boolean matrix $A$ is said to be 
\begin{enumerate}
\item convergent if and only if there exists in the sequence $\lbrace A^{k} : k \geq 1\rbrace$  a power $A^{m}$ of $A$ such that $A^{m}=A^{m+1}$.
\item oscillatory or periodic if and only if there exists in the sequence $\lbrace A^{k} : k \geq 1\rbrace$, a power $A^{m}$ of $A$ such that $A^{m}=A^{m+q}$ where $q$ is the smallest integer for which this holds and $q > 1$.
\end{enumerate}
\end{definition}

\begin{theorem} 
\label{Rosenblatt-1}
If $A$ is a Boolean matrix of order $n$ with graph $\G(A)$, then the following hold:
\begin{enumerate}
\item $A^{k}$ converges to the zero matrix if and only if $\G(A)$ contains no cyclic nets.
\item $A^{k}$ converges to $J$ (the matrix of all $1$'s) if and only if $\G(A)$ is a universal cyclic net.
\item $A$ is oscillatory if and only if $\G(A)$ contains at least one maximal cyclic net which is not universal. The period of a submatrix corresponding to a specified non-universal maximal cyclic net in $\G(A)$ is given by the greatest common divisor of the order of all simple cyclic nets contained in the maximal net.
\end{enumerate} 
\end{theorem}

\begin{remark}
Cyclic nets are the connected components of $\G(A)$ and the order of a cyclic net is the number of vertices in the cyclic net. A cyclic net is said to be a universal cyclic net if for some positive integer $q$ every point of the cyclic net is attainable in $q$ steps from some fixed point in the cyclic net.   
\end{remark} 

\noindent
We now define asymptotic period of a sequence of matrices, as given in \cite{Pang-Guu}.

\begin{definition} 
\label{fourpointtwo}
A sequence $\lbrace A_{k}:\ k \in \mathbb{N}\rbrace$ of matrices in $M_{n} (\reals)$ is called asymptotic $q$-periodic if 
\begin{equation} 
\label{Asy Period}
\lim\limits_{k\, \to\, \infty} A_{j + k q}\ \ =\ \ \widetilde{A^{(j)}}
\end{equation}
exists for $j= 1, 2, \ldots, q$. The minimal such $q$ is called the asymptotic period of the sequence.  
\end{definition}

\noindent
We note that the Definitions \eqref{fourpointone} and \eqref{fourpointtwo} coincide, for Boolean matrices. 
\medskip 

\noindent 
Let $J_{n}$ be the square matrix of ones of appropriate size and $A$ be an $n\times n$ nonnegative matrix with $A \leq J_{n}$ (entrywise). Consider the following decomposition:
\begin{equation} 
\label{Decomposition}
A\ \ =\ \ B(A) \ \oplus \ R(A),
\end{equation}
where
\begin{tabular}{ c r c l}
$B(A)$ = $ \begin{cases}
					 a_{ij} \quad &\text{if}\  \, a_{ij}=1  \\ 
					 0 \quad &\text{if}\ \, a_{ij}<1 \\ 
				 \end{cases} $
		  & and  & 
$R(A)$ = $\begin{cases} 
					0 \quad &\text{if}\ \, a_{ij}=1  \\ 
					a_{ij} \quad &\text{if}\ \, a_{ij}<1 \\ 
				 \end{cases} $ 
\end{tabular}.
\medskip 

\noindent 
The following result is due to Pang and Guu. 

\begin{theorem} \cite{Pang-Guu} 
\label{PG-result}
Let $A=B(A) \oplus R(A)$. $B(A)$ has period $q$ if and only if the sequence $\lbrace A^{k}\ :\ k \in \mathbb{N} \rbrace$ has asymptotic period $q$.
\end{theorem}

\noindent
\begin{remark} 
\label{Boolean Period}
From the above definitions, results and remarks, it follows that the period of a Boolean matrix is the same as the greatest common divisor of all possible lengths of simple circuits in $\G(A^{C})$.  
\end{remark}

\section{Proof of Theorem \eqref{result 1}}

In this section, we prove Theorem \eqref{result 1}. We begin with an extra assumption that $A$ is irreducible. The reducible case will follow after the proof of the irreducible case. 
\medskip 

\noindent 
\begin{proof}[of Theorem \eqref{result 1}] Let $A$ be a given irreducible matrix with $\mu(A) \le 1$. By Lemma \eqref{DAD}, without loss of generality, we can assume $A \leq J_{n}$. By equation \eqref{Decomposition}, we have $A=B(A) \oplus R(A)$. Since $B(A)$ is a Boolean matrix, there exists $q \ge 1$ such that $B(A)$ is $q$-periodic. Applying Theroem \eqref{PG-result} to $A$, the sequence $\lbrace A^{k} : \ k\in\mathbb{N}\rbrace$ has asymptotic period $q$. That is, for each $j = 1, 2, \ldots, q,\ \lim\limits_{k\, \to\, \infty} A^{j + k q} = \widetilde{A^{(j)}}$, exists. 
\medskip 

\noindent 
For $x \in \R$, define $\xi_{x}^{(j)} := \widetilde{A^{(j)}} \otimes x$. Now we show that $\xi_{x}^{(j)}$ is a periodic point of $A$ with its period dividing $q$. Consider 
\[ A^{q} \otimes \xi_{x}^{(j)}\ =\ A^{q} \otimes \lim\limits_{k\, \to\, \infty} A^{k q + j} \otimes x\ =\ \lim\limits_{k\, \to\, \infty} A^{k q + q + j} \otimes x\ =\ \lim\limits_{k\, \to\, \infty} A^{(k + 1) q + j} \otimes x\ =\ \xi_{x}^{(j)}. \]
Hence $\xi_{x}^{(j)}$ is a periodic point of $A$ with its period dividing $q$.
\medskip

\noindent
We now continue with $A$ being reducible. Then there exists a permutation matrix $P$ such that $PAP^{T}$ is a block triangular matrix in its Fr\"{o}benius normal form as given in \eqref{Frobenius normal form}. 
\begin{equation}
PAP^{T}\ \ = \ \  \left[\begin{matrix}
						A_{11}	&	A_{12}	& \cdots	& A_{1m}\\
						0		&	A_{22}	& \cdots	& A_{2m}\\
						\vdots	&	\vdots	& \ddots	& \vdots\\
						0		&	0		& \cdots	& A_{mm}\\
 					\end{matrix}\right],
\end{equation}
where each block $A_{ii}$ is square and is either irreducible or a $1 \times 1$ null matrix.
\medskip

\noindent
Suppose $A_{ii}$ is an irreducible matrix. Then, applying Theorem \eqref{DAD} to $A_{ii}$, there exists a diagonal matrix $D_{i}$ such that $D_{i}^{-1}A_{ii}D_{i} \leq J_{n_{i}}$, where $n_{i}$ is order of $A_{ii}$. However, if $A_{ii}$ is a $1 \times 1$ null matrix, we choose $D_{i} = \left[ 1 \right]$. Now, define $D$ as the direct sum of $D_{1}, D_{2}, \cdots, D_{m}$. Then, $D^{-1}\D D \leq J$ which is defined as the direct sum of $J_{n_{i}}$ for $1 \le i \le m$. Without loss of generality, we assume $\D \leq  J \leq J_{n}$.  				
\medskip

\noindent
Define $\Lambda\ =\ \lbrace i \in \lbrace 1, 2, \ldots, m \rbrace\ :\ \mu(A_{ii}) = 1 \rbrace$. Then, for $i \in \Lambda$, we apply Theorem \eqref{Elsner-result} to deduce that there exist $q^{(i)}$ and $ t_{i}$ such that $A_{ii}^{t + q^{(i)}} = A_{ii}^{t}$ for all $t \geq t_{i}$. However, for $i \notin \Lambda$, given $\epsilon > 0$, there exists an integer $t_{i}$ such that $A_{ii}^{t} < \epsilon  J_{n_{i}}$ for all $t\geq t_{i}$. Thus, $\lim\limits_{k\, \to\, \infty} A_{ii}^{k} = 0$ that yields $q^{(i)} = 1$, in this case. Now, defining $q = {\rm lcm} (q^{(1)}, q^{(2)}, \ldots, q^{(m)})$ and choosing $t_{0}$ to be bigger than $\max \{ t_{1}, t_{2}, \ldots, t_{m} \}$, we obtain for all $t \ge t_{0}$, 
\begin{equation} 
\label{Aii period} 
A_{ii}^{t + q} = A_{ii}^{t}\ \ \ \text{for}\ i \in \Lambda\ \ \ \text{and}\ \ \ A_{ii}^{t} < \epsilon J_{n_{i}}\ \ \ \text{for}\ i \notin \Lambda. 
\end{equation} 

\noindent 
Call $\lim\limits_{k\, \to\, \infty} A_{ii}^{kq} = \widetilde{A_{i}}$ for $1 \le i \le m$. For any vector $x^{(n_{i})} \in \mathbb{R}^{n_{i}}_{+}$, let $\lim\limits_{k\, \to\, \infty} A_{ii}^{kq} \otimes x^{(n_{i})} = \xi_{i}$. As $\D$ is the direct sum of $A_{ii}$, we obtain $\lim\limits_{k\, \to\, \infty} \D^{k q}$ to be the direct sum of the matrices $\widetilde{A_{i}}$ and for any vector $x = \left( x^{(n_{1})}, \ldots, x^{(n_{m})} \right) \in \R$, 
\begin{equation}
\lim\limits_{k\, \to\, \infty} \D^{k q} \otimes x\ \ =\ \ \xi_{x},\ \ \ \ \ \ \text{where}\ \ \xi_{x}\ = \left( \xi_{1}, \xi_{2}, \ldots, \xi_{m} \right). 
\end{equation}	

\noindent
Owing to the expression $A = \D \oplus \N$, we have for any $k \in \mathbb{N}$, 
\[ A^{k}\ \ \leq\ \ \D^{k} \oplus \N  \oplus \N^{2} \oplus \cdots \oplus \N^{k}. \] 
Since the index of nilpotency of $\N$ is at most $n$, we have for any $k \ge n$, 
\begin{equation} 
\label{Bound}
A^{k}\ \ \leq \ \ \D^{k}\, \oplus\, \N\, \oplus\, \N^{2}\, \oplus\, \cdots\, \oplus\, \N^{n - 1}. 
\end{equation}
Hence 
\[ \lim\limits_{k\, \to\, \infty} \D^{kq}\ \leq\ \lim\limits_{k\, \to\, \infty} A^{kq}\ \leq\ \lim\limits_{k\, \to\, \infty} \left[ \D^{kq}\, \oplus\, \N\, \oplus\, \N^{2}\, \oplus\, \cdots\, \oplus\, \N^{n - 1} \right]. \]

\noindent
We now prove the following claim. 
\[ \lim\limits_{k\, \to\, \infty} A^{k q}\ \text{exists if and only if} \lim\limits_{k\, \to\, \infty} \D^{k q}\ \text{exists}. \] 
Suppose $\lim\limits_{k\, \to\, \infty} A^{k q}$ exists, then it is trivially true that $\lim\limits_{k\, \to\, \infty} \D^{k q}$ exists. Thus, it is sufficient to prove that $\lim\limits_{k\, \to\, \infty} A^{k q}$ exists, when $\lim\limits_{k\, \to\, \infty} \D^{k q}$ exists. We prove this by induction on $m$. 
\medskip 

\noindent
When $m = 2$ with \begin{tabular}{ c r c l} $\D = \left[ \begin{matrix} A_{11} & 0 \\ 0 & A_{22} \end{matrix} \right]$ and $\N = \left[ \begin{matrix} 0 & A_{12} \\ 0 & 0 \end{matrix} \right]$, \end{tabular} where $A_{11}$ and $A_{22}$ are irreducible matrices with entries less than or equal to $1$ and $\N$ is a nilpotent matrix of nilpotency index $2$, we have for $k \ge 2$, 
\begin{eqnarray*}
A^{k}	&	=	&	\left[ \begin{matrix}
					A_{11}^{k} & 0 \\
					0	& A_{22}^{k}	
					\end{matrix} \right] \oplus \left[ \begin{matrix}
					0 	& \bigoplus\limits_{\ell\, =\, 0}^{k - 1} A_{11}^{\ell} A_{12} A_{22}^{k - \ell - 1} \\
 				0	& 0	
					\end{matrix} \right].
\end{eqnarray*}

\subsubsection*{Case 1: Suppose $\mu(A_{11}) = \mu(A_{22}) = 1$.} 

Consider
\begin{eqnarray*} 
\bigoplus_{\ell\, =\, 0}^{kq - 1} A_{11}^{\ell} A_{12} A_{22}^{kq - \ell - 1} & = & \bigoplus_{\ell\, =\, 0}^{t_{0} - 1} A_{11}^{\ell} A_{12} A_{22}^{kq - \ell - 1}\; \bigoplus_{\ell\, =\, t_{0}}^{kq - t_{0} - 1} A_{11}^{\ell} A_{12} A_{22}^{kq - \ell - 1}\; \bigoplus_{\ell\, =\, kq - t_{0}}^{kq - 1} A_{11}^{\ell} A_{12} A_{22}^{kq - \ell - 1} \\ 
& =: & I_{1}(k)\,\oplus I_{2}(k)\,\oplus I_{3}(k), 
\end{eqnarray*} 
where $q$ and $t_{0}$ are determined as in Equation \eqref{Aii period}. For $kq \geq t_{0}$ and for $1 \leq \ell \leq q$, we observe that $A_{ii}^{(k + 1)q - \ell} = A_{ii}^{t_{0} + q -\ell}$. Note that $I_{j} (k)$ is independent of $k$; the arguments for $I_{1}$ and $I_{2}$ are presented below. 
\begin{eqnarray*}
I_{1}(k) & = & \left( A_{12} A_{22}^{kq - 1} \right) \oplus \left( A_{11} A_{12} A_{22}^{kq - 2} \right)\; \oplus\; \cdots\; \oplus\; \left( A_{11}^{t_{0} - 1} A_{12} A_{22}^{kq - t_{0}} \right) \\ 
& = & \left( A_{12} A_{22}^{t_{0} + q - 1} \right) \oplus \left( A_{11} A_{12} A_{22}^{t_{0} + q - 2} \right)\; \oplus\; \cdots\; \oplus\; \left( A_{11}^{t_{0} - 1} A_{12} A_{22}^{t_{0}} \right), 
\end{eqnarray*}
whereas, 
\begin{eqnarray*}
I_{2}(k) & = & \left( A_{11}^{t_{0}} A_{12} A_{22}^{kq - t_{0} - 1} \right)\; \oplus\; \left( A_{11}^{t_{0} + 1} A_{12} A_{22}^{kq - t_{0} - 2} \right)\; \oplus\; \cdots\; \oplus\; \left( A_{11}^{kq - t_{0} - 1} A_{12} A_{22}^{kq - t_{0}} \right) \\
& = & \left( A_{11}^{t_{0}} A_{12} A_{22}^{t_{0} + q - 1} \right)\; \oplus\; \left( A_{11}^{t_{0} + 1} A_{12} A_{22}^{t_{0} + q - 2} \right)\; \oplus\; \cdots\; \oplus\; \left( A_{11}^{t_{0} + q - 1} A_{12} A_{22}^{t_{0}} \right). 
\end{eqnarray*}

\noindent 
Thus, for $k$ sufficiently large, the nilpotent part of the sequence $\lbrace A^{kq}\rbrace_{k\in\mathbb{Z}_{+}}$ is independent of $k$. 

\subsubsection*{Case 2: Suppose $\mu(A_{11}) = 1$ while $\mu(A_{22}) < 1$.}

In this case, we break the sum as 
\begin{eqnarray*}
\bigoplus_{\ell\, =\, 0}^{kq - 1} A_{11}^{\ell} A_{12} A_{22}^{kq - \ell - 1} & = & \bigoplus_{\ell\, =\, 0}^{(k - 1)q - t_{0} - 1} A_{11}^{\ell} A_{12} A_{22}^{kq - \ell - 1}\; \bigoplus_{\ell\, =\, (k - 1)q - t_{0}}^{kq - t_{0} - 1} A_{11}^{\ell} A_{12} A_{22}^{kq - \ell - 1}\; \bigoplus_{\ell\, =\, kq - t_{0}}^{kq - 1} A_{11}^{\ell} A_{12} A_{22}^{kq - \ell - 1} \\
& =: & I_{4}(k)\, \oplus I_{5}(k)\, \oplus I_{6}(k), 
\end{eqnarray*}
where $q$ and $t_{0}$ are determined as in Equation \eqref{Aii period}. Defining $\alpha_{ij} := \left( I_{5}(k) \oplus I_{6} (k) \right)_{ij},\ \beta := \max\limits_{i, j} \left( A_{12} \right)_{ij}$ and choosing $\epsilon < \min\limits_{i, j} \left\{ \alpha_{ij} \beta^{-1} \right\}$, we observe that $\left( I_{4} (k) \right)_{ij} \le \left( I_{5}(k) \oplus I_{6} (k) \right)_{ij}$, thus proving that $I_{4}$ does not contribute to the sum above. 
\medskip 

\noindent 
Using Equation \eqref{Aii period} for the irreducible matrix $A_{11}$ whose circuit geometric mean is $1$, we see that 
\begin{eqnarray*} 
I_{5} (k) & = & \left( A_{11}^{(k - 1) q - t_{0}} A_{12} A_{22}^{t_{0} + q - 1} \right)\; \oplus\; \left( A_{11}^{(k - 1) q - t_{0} + 1} A_{12} A_{22}^{t_{0} + q - 2} \right)\; \oplus\; \cdots\; \oplus\; \left( A_{11}^{kq - t_{0} - 1} A_{12} A_{22}^{t_{0}} \right) \\ 
& = & \left( A_{11}^{t_{0} + r_{1}} A_{12} A_{22}^{t_{0} + q - 1} \right)\; \oplus\; \left( A_{11}^{t_{0} + r_{2}} A_{12} A_{22}^{t_{0} + q - 2} \right)\; \oplus\; \cdots\; \oplus\; \left( A_{11}^{t_{0} + r_{q}} A_{12} A_{22}^{t_{0}} \right), 
\end{eqnarray*} 
whereas, 
\begin{eqnarray*} 
I_{6} (k) & = & \left( A_{11}^{kq - t_{0}} A_{12} A_{22}^{t_{0} - 1} \right)\; \oplus\; \left( A_{11}^{kq - t_{0} + 1} A_{12} A_{22}^{t_{0} - 2} \right)\; \oplus\; \cdots\; \oplus\; \left( A_{11}^{kq - 1} A_{12} \right) \\ 
& = & \left( A_{11}^{t_{0} + s_{1}} A_{12} A_{22}^{t_{0} - 1} \right)\; \oplus\; \left( A_{11}^{t_{0} + s_{2}} A_{12} A_{22}^{t_{0} - 2} \right)\; \oplus\; \cdots\; \oplus\; \left( A_{11}^{t_{0} + s_{t_{0}}} A_{12} \right), 
\end{eqnarray*} 
where $0 \le r_{i}, s_{j} \le q - 1$, thus making both $I_{5}$ and $I_{6}$ independent of $k$, for $k$ sufficiently large. 
\medskip 

\noindent 
Since the cases discussed above exhaust the possibilities for the representation of the matrix $A$ in its Fr\"{o}benius normal form, we conclude by virtue of our hypothesis, that $\lim\limits_{k\, \to\, \infty} A^{kq} = \widetilde{A} = \left[ \begin{matrix} A_{11}^{t_{0}} & M \\ 0 & A_{22}^{t_{0}} \end{matrix} \right]$ exists. We now prove that there exists some $t' > 0$ such that $\left( \widetilde{A} \right)^{t + q}  = \left( \widetilde{A} \right)^{t}$ for all $t > t'$, thus ensuring that the period of $\widetilde{A}$ divides $q$. 
\medskip 

\noindent 
Readers may find that the arguments presented to prove periodicity are analogous to the ones presented earlier to prove the existence of the limit; {\it i.e.}, a clever usage of Elsner's Theorem \eqref{Elsner-result} and a certain splitting of the max sum and comparison arguments. In the case when $\mu(A_{11}) = 1$ and $\mu(A_{22}) < 1$, we have 
\[ \left( \widetilde{A} \right)^{t + q}\ \ =\ \ \left[ \begin{matrix} A_{11}^{t_{0} (t + q)} & A_{11}^{t_{0} (t + q - 1)} M \\ 0 & 0 \end{matrix} \right]\ \ =\ \ \left[ \begin{matrix} A_{11}^{t_{0} t} & A_{11}^{t_{0} (t - 1)} M \\ 0 & 0 \end{matrix} \right]\ \ =\ \ \left( \widetilde{A} \right)^{t}\ \ \forall t > t' = 1. \] 
In the case when $\mu(A_{11}) = \mu(A_{22}) = 1$, we have 
\begin{eqnarray*} 
\left( \widetilde{A} \right)^{t + q} & = & \left[ \begin{matrix} A_{11}^{t_{0} (t + q)} & \bigoplus\limits_{l\, =\, 0}^{t + q - 1} A_{11}^{t_{0} l} M A_{22}^{t_{0} (t + q - l - 1)} \\ 0 & A_{22}^{t_{0} (t + q)} \end{matrix} \right] \\ 
& = & \left[ \begin{matrix} A_{11}^{t_{0} t} & 0 \\ 0 & A_{22}^{t_{0} t} \end{matrix} \right] \oplus \left[ \begin{matrix} 0 & \bigoplus\limits_{l\, =\, 0}^{t - 1} A_{11}^{t_{0} l} M A_{22}^{t_{0} (t + q - l - 1)}\ \bigoplus\limits_{l\, =\, t}^{t + q - 1} A_{11}^{t_{0} l} M A_{22}^{t_{0} (t + q - l - 1)} \\ 0 & 0 \end{matrix} \right] \\ 
& = & \left[ \begin{matrix} A_{11}^{t_{0} t} & 0 \\ 0 & A_{22}^{t_{0} t} \end{matrix} \right] \oplus \left[ \begin{matrix} 0 & \bigoplus\limits_{l\, =\, 0}^{t - 1} A_{11}^{t_{0} l} M A_{22}^{t_{0} (t - l - 1)} \\ 0 & 0 \end{matrix} \right] \\ 
& = & \left( \widetilde{A} \right)^{t}\ \ \ \ \forall t > t'\ge q. 
\end{eqnarray*} 

\noindent 
It is an easy exercise to note that one may analogously prove the above statements even when we consider $A^{kq + j}$ whose limit, we denote by $\widetilde{A^{(j)}}$, as $k \to \infty$. Moreover, for any $x \in \R$, 
\[ \lim\limits_{k\, \to\, \infty} A^{kq + j} \otimes x\ \ =\ \ \lim\limits_{k\, \to\, \infty} \left[ \begin{matrix} \left( A_{11}^{kq + j} \otimes x_{1} \right) \oplus \left( \bigoplus\limits_{\ell\, =\, 0}^{kq + j - 1} A_{11}^{\ell} A_{12} A_{22}^{kq + j - \ell - 1} \otimes x_{2} \right) \\ A_{22}^{kq + j} \otimes x_{2} \end{matrix} \right]\ \ =\ \ \left[ \begin{matrix} \xi_{1}^{(j)} \oplus \eta_{1}^{(j)} \\ 0 \end{matrix} \right], \] 
thereby, completing the proof in the base case when $m = 2$. A simple induction argument on $m$, now clinches the proof. 
\end{proof} 

\section{Proofs of Theorems \eqref{result 2} and \eqref{result 3}}  

We start this section with the hypothesis that $\C$ is a collection of pairwise commuting nonnegative matrices with $\mu(A_{i}) \leq 1$. Consider a $p$-lettered word $\omega \in \Sigma_{N}^{p}$, possibly with all letters present at least once. Let $A_{\omega}$ be the matrix associated with the $p$-lettered word $\omega$, as defined in Equation \eqref{matrixproduct}. From Theorem \eqref{commuting matrices mu inequality}, we know that 
\[ \mu(A_{\omega})\ \ \leq\ \ \mu(A_{\omega_{p}}) \times \mu(A_{\omega_{p - 1}}) \times \ldots \times \mu(A_{\omega_{1}})\ \ \leq\ \ 1. \] 

\begin{proof} [of Theorem \eqref{result 2}] 
Applying Theorem \eqref{result 1} to $A_{i}$, we obtain an integer $q_{i}$ such that for every $1 \le j \le q_{i}$, we have
\[\lim\limits_{k\, \to\, \infty} A_{i}^{k q_{i} + j}\ \ =\ \ \widetilde{A^{(j)}_{i}}. \] 

\noindent
Define $q := {\rm lcm} \left( q_{1},\, q_{2},\, \cdots,\, q_{N} \right)$. As the matrices commute, we see from the associativity of $\otimes$ that 
\[A_{\omega}\ \ =\ \ A_{1}^{p_{1}} \otimes A_{2}^{p_{2}} \otimes \cdots \otimes A_{N}^{p_{N}},\]
where $p_{i}$ is the number of occurrences of $i$ in the $p$-lettered word $\omega$ with $p = p_{1} + p_{2} + \cdots + p_{N}$. Observe that $A_{\omega}^{k} = A_{1}^{kp_{1}} \otimes A_{2}^{kp_{2}} \otimes \cdots \otimes A_{N}^{kp_{N}}$ from which, it follows that 
\[\lim\limits_{k\, \to \infty} A_{\omega}^{kq + j}\ \ =\ \ \bigotimes\limits_{i\, =\, 1}^{N} \left( \widetilde{A^{(j)}_{i}} \right)^{p_{i}}. \] 
\noindent 
Finally, by our definition of $q$, we observe that it is independent of the choice of the word $\omega$ as well as its length. For any word $\omega$ that may not contain all the letters, one makes an analogous argument with a set that has fewer letters. 
\end{proof}
\medskip 

\noindent
The following is an obvious corollary to the above theorem. 

\begin{corollary}
\label{coro1}
Let $A_{1}$ and $A_{2}$ be two commuting matrices with all the max eigenvalues being either $0$ or $1$. Suppose $\omega \in \Sigma_{2}^{p}$ is such that the letter $i$ occurs $p_{i}$ many times and has the associated matrix product $A_{\omega}$. Then, there exist $q$ and $t_{0}$ such that for any $1 \le j \le q$, we have 
\[ \lim_{k\, \to\, \infty} A_{\omega}^{kq + j}\ \ =\ \ L_{\omega}^{(j)}, \] 
where $L_{\omega}^{(j)} \in \left\{ A_{1}^{p_{1} t_{0}} \otimes A_{2}^{p_{2} t_{0}},\, A_{1}^{p_{1}(t_{0} + 1)} \otimes A_{2}^{p_{2}(t_{0} + 1)},\, \cdots,\, A_{1}^{p_{1}(t_{0} + q - 1)} \otimes A_{2}^{p_{2}(t_{0} + q - 1)}\right\}$.  
\end{corollary}

\noindent
We now prove Theorem \eqref{result 3} with the supposition that the collection $\C$ is not necessarily pairwise commuting. A key hypothesis in the proof of Theorem \eqref{result 3} is the existence of a set of common max eigenvectors for the collection $\C$. 
\medskip
 
\noindent 
\begin{proof} [of Theorem \eqref{result 3}]
Let $E = \lbrace v_{1},\, v_{2},\, \cdots,\, v_{m} \rbrace$ be the set of common max eigenvectors for the collection $\C$. We rearrange the eigenvectors (and rename them, if necessary) as 
\[ \left\{ v_{1},\, v_{2},\, \cdots,\, v_{\kappa},\, v_{\kappa + 1},\, v_{\kappa + 2},\, \cdots,\, v_{m}\right\} \] 
so that $A_{i} \otimes v_{j} = v_{j}$ for all $1 \leq i \leq N$ and $1 \leq j \leq \kappa$. Assuming that the word $\omega \in \Sigma_{N}^{p}$ contains all the letters at least once, it follows that for every $\kappa + 1 \leq j \leq m,\ \lim\limits_{k\, \to\, \infty} A_{\omega}^{k} \otimes v_{j} = 0$. Now, for any $x = \alpha_{1} v_{1} \oplus \alpha_{2} v_{2} \oplus \cdots \oplus \alpha_{m} v_{m} \in \mathcal{LC} (E)$, consider
\begin{eqnarray*}
\lim\limits_{k\, \to\, \infty} A_{\omega}^k\otimes x & = & \left[ \alpha_{1} \lim\limits_{k\, \to\, \infty}  \left( A_{\omega}^{k} \otimes v_{1} \right) \right] \oplus \left[ \alpha_{2} \lim\limits_{k\, \to\, \infty} \left( A_{\omega}^{k} \otimes v_{2} \right) \right] \oplus \cdots \oplus \left[ \alpha_{m} \lim\limits_{k\, \to\, \infty} \left( A_{\omega}^{k} \otimes v_{m} \right) \right] \\ 
& = & \alpha_{1} v_{1} \oplus \alpha_{2} v_{2} \oplus \cdots \oplus \alpha_{\kappa} v_{\kappa} \\ 
& =: & \xi_{x}, 
\end{eqnarray*}  
with $\xi_{x}$ being a fixed point for every matrix in the collection $\C$.
\end{proof}

\section{Examples}

In this section, we illustrate the theorems in this article with a few examples. We begin with an example that concerns Theorem \eqref{result 1}.

\begin{example}
\label{example1}
Let $A\ =\ \left[ \begin{matrix} 0.2 & 1  & 4 \\ 1  & 0.5 & 6 \\ 0  & 0  & 0.9 \end{matrix} \right]$ with corresponding graph given by 
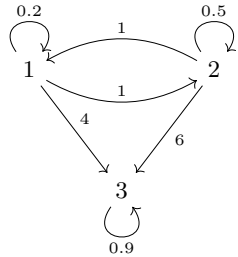
\begin{center}
\small
\begin{tikzcd}
1 \arrow["0.2", loop, distance=2em, in=55, out=125] \arrow[rdd, "4"] \arrow[rr, "1", bend right] &                                                       & 2 \arrow[ll, "1"', bend right] \arrow[ldd, "6"] \arrow["0.5", loop, distance=2em, in=55, out=125] \\
                                                                                                 &                                                       &                                                                                                   \\
                                                                                                 & 3 \arrow["0.9"', loop, distance=2em, in=305, out=235] &                                                                                                  
\end{tikzcd}
\end{center}
Recall that the critical circuit is the one where the maximum circuit geometric mean is attained in $\G(A)$; in this case, it is $(2, 1) (1, 2)$, that yields $\mu (A) = 1$ and using Remark \eqref{Boolean Period}, we also obtain $q = 2$. Then, for any $k \in \mathbb{Z}_{+}$, we have 
\[ A^{2k + 1}\ =\ \left[ \begin{matrix} 0.5 & 1 & 5.4 \\ 1  & 0.5 & 6 \\ 0  & 0  & (0.9)^{2k + 1} \end{matrix} \right]\ \ \ \text{while}\ \ \ A^{2k}\ =\ \left[ \begin{matrix} 1 & 0.5  & 6 \\ 0.5  & 1 & 5.4 \\ 0  & 0  & (0.9)^{2k} \end{matrix} \right]. \] 
Thus, 
\[ \widetilde{A^{(1)}}\ =\ \left[ \begin{matrix} 0.5 & 1 & 5.4 \\ 1  & 0.5 & 6 \\ 0  & 0  & 0 \end{matrix} \right]\ \ \ \text{while}\ \ \ \widetilde{A^{(2)}}\ =\ \left[ \begin{matrix} 1 & 0.5  & 6 \\ 0.5  & 1 & 5.4 \\ 0  & 0  & 0 \end{matrix} \right]. \] 
Hence, for any $x = (x_{1}, x_{2}, x_{3})^{t} \in \mathbb{R}^{3}_{+}$, it is easy to verify that 
\[ \lim_{k\, \to\, \infty} A^{2k + 1} \otimes x\ \ =\ \ \begin{pmatrix} 0.5 x_{1} \oplus x_{2} \oplus 5.4 x_{3} \\ x_{1} \oplus 0.5 x_{2} \oplus 6 x_{3} \\ 0 \end{pmatrix} \ \ \ \text{whereas}\ \ \ \lim_{k\, \to\, \infty} A^{2k} \otimes x\ \ =\ \ \begin{pmatrix} x_{1} \oplus 0.5 x_{2} \oplus 6 x_{3} \\ 0.5 x_{1} \oplus x_{2} \oplus 5.4 x_{3} \\ 0 \end{pmatrix}. \]
\end{example}
\medskip

\noindent
The following example illustrates Theorem \eqref{result 2} and Corollary \eqref{coro1}. 

\begin{example}
\label{example2}
Consider the collection $\{ A_{1}, A_{2}, A_{3} \}$ of pairwise commuting $5 \times 5$ nonnegative matrices given by 
\[ A_{1}\ \ =\ \ \left[\begin{matrix} 0 & 1 & 0 & 8 & 5 \\ 0	& 0 & 1 & 5 & 8 \\ 1 & 0 & 0 & 8 & 5 \\ 0 & 0 & 0 & 1 & 0.5 \\ 0	& 0 & 0 & 0.5 & 1 \end{matrix}\right],\ \ \ \  A_{2}\ \ =\ \ \left[\begin{matrix}0 & 0 & 1 & 9 & 9 \\ 1 & 0 & 0 & 9 & 9 \\ 0 & 1 & 0 & 9 & 9 \\ 0 & 0 & 0 & 1 & 1 \\ 0 & 0 & 0 & 1 & 1 \end{matrix}\right]\ \ \ \text{and}\ \ \ \ A_{3}\ \ =\ \ \left[\begin{matrix}1 & 0 & 0 & 8 & 9 \\ 0 & 1 & 0 & 8 & 9 \\ 0 & 0 & 1 & 8 & 9 \\ 0 & 0 & 0 & 0.8 & 1 \\ 0 & 0 & 0 & 1 & 0.8 \end{matrix}\right] \]
with $\mu(A_{i})\; =\; 1$ for all $i$. Observe that all the max eigenvalues of the three matrices are equal to $1$. Then, owing to the Theorem \eqref{result 1}, we find that the matrices $A_{1},\, A_{2}$ and $A_{3}$ have asymptotic periods $q_{1},\, q_{2}$ and $q_{3}$ given by $3,\, 3$ and $2$ respectively, {\it i.e.}, $\lim\limits_{k\, \to\, \infty} A_{i}^{k q_{i} + j} = \widetilde{A_{i}^{(j)}}$.
\medskip 

\noindent
We now discuss all possible cases of $\omega \in \Sigma_{3}^{2}$, with distinct letters. 
\begin{enumerate} 
\item[(i)] Let $A_{\omega} = A_{1} \otimes A_{2}$. In this case, we find that $\lim\limits_{k\, \to\, \infty} A_{\omega}^{k q_{\omega} + j}\ =\ A_{\omega}$ (irrespective of $j$) is a matrix with period $1$, while $q_{\omega} = 3$. 

\item[(ii)] Let $A_{\omega} = A_{2} \otimes A_{3}$. Then, $\lim\limits_{k\, \to\, \infty} A_{\omega}^{kq_{\omega} + j}\ =\ L_{\omega}^{(j)}$, where $L_{\omega}^{(j)}$ has period $3$ while $q_{\omega} = 6$. Further, $L_{\omega}^{(j)}\ \in\ \left\lbrace A_{2} \otimes A_{3},\, A_{2}^{2} \otimes A_{3}^{2},\, A_{2}^{3} \otimes A_{3}^{3} \right\rbrace$. 

\item[(iii)] Let $A_{\omega} = A_{1} \otimes A_{3}$. In this case, $\lim\limits_{k\, \to\, \infty} A_{\omega}^{kq_{\omega} + j}\ =\ L_{\omega}^{(j)}$, where $L_{\omega}^{(j)}$ has period $q_{\omega} = 6$. Further, $L_{\omega}^{(j)}\ \in\ \left\lbrace A_{1} \otimes A_{3},\, A_{1}^{2} \otimes A_{3}^{2},\, \cdots,\, A_{1}^{5} \otimes A_{3}^{5} \right\rbrace$. 
\end{enumerate} 
\medskip 

\noindent
For a general $p$-lettered word $\omega$ that contains each of the letters at least once, suppose $p_{i}$ counts the number of occurrences of the letter $i$ in $\omega$, we find that $\lim\limits_{k\, \to\, \infty} A_{\omega}^{kq + j}= L_{\omega}^{(j)}$. One obtains the period of $L_{\omega}^{(j)}$ to be either $1$ or $3$ while $q = 6$, here. This is a special occasion since the second diagonal block of the matrices $A_{2}$ and $A_{3}$, expressed in their Fr\"{o}benius normal form, when multiplied yields $J_{2}$. Thus, the period of $A_{3}$ makes no contribution in the computation of the period of $L_{\omega}^{(j)}$. Hence, 
\[ L_{\omega}^{(j)}\ \ \in\ \ \left\lbrace A_{1}^{p_{1}} \otimes A_{2}^{p_{2}} \otimes A_{3}^{p_{3}},\ A_{1}^{2 p_{1}} \otimes A_{2}^{2 p_{2}} \otimes A_{3}^{2 p_{3}},\ A_{1}^{3 p_{1}} \otimes A_{2}^{3 p_{2}} \otimes A_{3}^{3 p_{3}} \right\rbrace. \]    
\end{example}
\medskip

\noindent
We conclude this section and the article with an example in the non-commuting case, that explains Theorem \eqref{result 3}.

\begin{example}
\label{example4}
Consider the matrices $A_{1}$ and $A_{2}$, where 
\[ A_{1}\ \ =\ \ \left[\begin{matrix} 0.9 & 0.45 & 5 & 6 & 27 \\ 0.45 & 0.9 & 1 & 23 & 8 \\ 0 & 0 & 0.9 & 1 & 0 \\ 0 & 0 & 0 & 0.2 & 1 \\ 0 & 0 & 1 & 0 & 0 \end{matrix}\right]\ \ \ \ \text{and}\ \ \ \ A_{2}\ \ =\ \ \left[\begin{matrix} 1  & 1 & 27 & 6 & 2 \\ 0.5 & 1 & 17 & 3 & 23 \\ 0 & 0 & 0.4 & 1 & 0 \\ 0 & 0 & 1 & 0.8 & 1 \\ 0 & 0 & 0.9 & 1 & 0.2 \end{matrix}\right].\]
The max eigenvalues for the matrix $A_{1}$ are $1$ and $0.9$, while those of $A_{2}$ are $1$ and $1$, with the corresponding common max eigenvectors for the matrices being $u = (27,\, 23,\, 1,\, 1,\, 1)^{T}$ and $v = (2,\, 1,\, 0,\, 0,\, 0)^{T}$ respectively. Suppose $\omega \in \Sigma_{2}^{p}$ with the presence of both the letters at least once. For $x \in \mathcal{LC} (E)$ with $x\ =\ \alpha u \oplus \beta v$ for some $\alpha,\, \beta \geq 0$, we have 
\[ \lim\limits_{k\, \to\, \infty} A_{\omega}^{k} \otimes x\ \ =\ \ \xi_{x}\ \ =\ \ \alpha u, \]
where $\xi_{x}$ is a common fixed point of $A_{1}$ and $A_{2}$.
\end{example}

\bibliographystyle{amsplain}

\begin{thebibliography}{99}

\bibitem{Bapat}
{\sc Bapat, R. B.}, ``A max version of the Perron-Fr\"{o}benius theorem", \emph{Linear Algebra Appl.}, {\bf 275/376} (1998), 3 - 18.

\bibitem{Bapat-Stanford-van}
{\sc Bapat, R. B., Stanford, D. P.} and {\sc van den Driessche, P.}, ``Pattern properties and spectral inequalities in max algebra", \emph{SIAM J. Matrix Anal. Appl.}, {\bf 16} (1995), 964 - 976.

\bibitem{BP-book}
{\sc Berman, A.} and {\sc Plemmons, R. J.}, ``Nonnegative Matrices in the Mathematical Sciences", {\bf 9}, \emph{Society for Industrial and Applied Mathematics (SIAM)}, Philadelphia, PA, (1994).

\bibitem{Butkovic-Book}
{\sc Butkovic, P.}, ``Max-linear Systems: Theory and Algorithms", \emph{Springer Monographs in Mathematics}, {\bf 27}, (2010).

\bibitem{Elsner-van}
{\sc Elsner, L.} and {\sc van den Driessche, P.}, ``On the power method in max algebra", \emph{Linear Algebra Appl.}, {\bf 302/303} (1999), 17 - 32.

\bibitem{Elsner-van 1}
{\sc Elsner, L.} and {\sc van den Driessche, P.}, ``Modifying the power method in max algebra", \emph{Linear Algebra Appl.}, {\bf 332/334} (2001), 3 - 13. 

\bibitem{SJ-YP-SS-1}
{\sc Jayaraman, S., Prajapaty, Y. P.} and {\sc Sridharan S.}, ``Dynamics of products of nonnegative matrices", {\tt arXiv:2010.05560}.

\bibitem{LN-Book}
{\sc Lemmens, B.} and {\sc Nussbaum, R. D.}, ``Nonlinear Perron-Frobenius Theory", \emph{Cambridge Tracts in Mathematics},  {\bf 189}, Cambridge University Press, (2012).

\bibitem{RKatz}
{\sc Katz, R. D., Schneider, H.} and {\sc Sergeev, S.}, ``On commuting matrices in max algebra and in classical nonnegative algebra", \emph{Linear Algebra Appl.}, {\bf 436} (2012), 276 - 292.

\bibitem{Pang-Guu}
{\sc Pang, C. T.} and {\sc Guu, S. M.}, ``A note on the sequence of consecutive powers of a nonnegative matrix in max algebra", \emph{Linear Algebra Appl.}, {\bf 330} (2001), 209 - 213.

\bibitem{Rosenblatt}
{\sc Rosenblatt, D.}, ``On the graphs and asymptotic forms of finite Boolean relation matrices and stochastic matrices", \emph{Naval Res. Logist. Quart.}, {\bf 4} (1957), 151 - 167.

\end{thebibliography}

\end{document}